\RequirePackage{ifpdf}
\ifpdf 
\documentclass[pdftex]{sigma}
\else
\documentclass{sigma}
\fi

\input epsf.sty

\newtheorem{question}[theorem]{{\bf Question}}

\def\BN{\mathbb N}
\def\BZ{\mathbb Z}

\def\BQ{\mathbb Q}
\def\BR{\mathbb R}
\def\BC{\mathbb C}

\def\D{\Delta}

\def\M{\mathcal M}

\def\a{\alpha}

\def\l{\lambda}
\def\Ga{\Gamma}
\def\S{\Sigma}

\def\s{\sigma}

\def\ga{\gamma}

\def\la{\langle}
\def\ra{\rangle}

\def\gl{\mathfrak{gl}}

\def\e{\epsilon}
\def\Ga{\Gamma}
\def\d{\delta}
\def\b{\beta}

\def\s{\sigma}

\def\pt{\partial}

\def\ti{\widetilde}

\def\longto{\longrightarrow}

\def\hb{\hbar}
\def\Li{\mathrm{Li}}

\def\coeff{\mathrm{coef\/f}}
\def\fg{\mathfrak{g}}

\def\deg{\mathrm{deg}}

\def\longto{\longrightarrow}
\def\calA{\mathcal A}

\def\SL{\mathrm{SL}}
\def\gl{\mathfrak{gl}}
\def\sl{\mathfrak{sl}}
\def\Om{\Omega}
\def\PSU{\mathrm{PSU}}

\def\re{{\rm e}}

\numberwithin{equation}{section}

\numberwithin{theorem}{section}
\numberwithin{proposition}{section}
\numberwithin{lemma}{section}
\numberwithin{conjecture}{section}
\numberwithin{corollary}{section}
\numberwithin{definition}{section}
\numberwithin{remark}{section}

\begin{document}

\allowdisplaybreaks

\renewcommand{\thefootnote}{$\star$}

\renewcommand{\PaperNumber}{080}

\FirstPageHeading

\ShortArticleName{Analyticity of the Free Energy of a Closed 3-Manifold}

\ArticleName{Analyticity of the Free Energy of a Closed 3-Manifold\footnote{This paper
is a contribution to the Special Issue on Deformation
Quantization. The full collection is available at
\href{http://www.emis.de/journals/SIGMA/Deformation_Quantization.html}{http://www.emis.de/journals/SIGMA/Deformation\_{}Quantization.html}}}

\Author{Stavros GAROUFALIDIS~$^\dag$, Thang T.Q. L\^E~$^\dag$ and Marcos MARI\~NO~$^\ddag$}

\AuthorNameForHeading{S. Garoufalidis, T.T.Q. L\^e and M. Mari\~no}

\Address{$^\dag$~School of Mathematics, Georgia Institute of Technology,
         Atlanta, GA 30332-0160, USA}

\EmailD{\href{mailto:stavros@math.gatech.edu}{stavros@math.gatech.edu}, \href{mailto:letu@math.gatech.edu}{letu@math.gatech.edu}}

\Address{$^\ddag$~Section de Math\'ematiques, Universit\'e de Gen\`eve, CH-1211 Gen\`eve 4, Switzerland}
\EmailD{\href{mailto:Marcos.Marino@unige.ch}{Marcos.Marino@unige.ch}}

\ArticleDates{Received September 15, 2008, in f\/inal form November 06,
2008; Published online November 15, 2008}

\Abstract{The free energy of a closed 3-manifold is a 2-parameter formal power series
which encodes the perturbative Chern--Simons invariant (also known
as the LMO invariant) of a~closed 3-manifold with gauge group $U(N)$
for arbitrary $N$. We prove that the free energy of an arbitrary closed
3-manifold is uniformly Gevrey-$1$. As a corollary, it follows that the genus
$g$ part of the free energy is convergent in a neighborhood of
zero, independent of the genus. Our results follow from
an estimate of the LMO invariant, in a particular gauge, and from recent
results of Bender--Gao--Richmond on the asymptotics of the number of rooted maps
for arbitrary genus. We illustrate our results with an explicit formula
for the free energy of a Lens space. In addition, using the Painlev\'e
dif\/ferential equation, we obtain an asymptotic expansion for the number of
cubic graphs to all orders, stengthening the results of  Bender--Gao--Richmond.}

\Keywords{Chern--Simons theory; perturbation theory;
gauge theory; free energy; planar limit; Gevrey series; LMO invariant;
weight systems;  ribbon graphs; cubic graphs; lens spaces; trilogarithm;
polylogarithm; Painlev\'e I; WKB; asymptotic expansions; transseries;
Riemann--Hilbert problem}

\Classification{57N10; 57M25}

\section{Introduction}
\subsection{The free energy of a closed 3-manifold}
\label{sub.3man}

The free energy $F_M(\tau,\hb)$ of a closed 3-manifold $M$ (def\/ined
in equation \eqref{eq.FM2}) is a 2-parameter
formal power series with rational coef\/f\/icients in two variables
$\tau$ and $\hb^{2}$
\begin{gather}
\label{eq.FM}
F_M(\tau,\hb)=\sum_{g=0}^\infty F_{M,g}(\tau) \hb^{2g-2} \in \hb^{-2}
\BQ[[\tau,\hb^{2}]].
\end{gather}
The variable $\hb$ plays the role of Planck's constant, and the variable
$\tau=N \hb$ is the product of Planck's constant with $N$,
the size of the gauge group $U(N)$. The free energy
encodes the perturbative Chern--Simons invariant of $M$ along the trivial
f\/lat connection (also known as the LMO invariant of \cite{LMO}) with
gauge group $U(N)$ for arbitrary $N$. $F_{M,g}(\tau)$ (resp. $F_{M,0}(\tau)$)
is often called the {\em genus} $g$-contribution (resp. {\em planar limit})
to the free energy. Perturbative Quantum Field Theory
for a gauge group of f\/ixed size $N$ typically leads to factorially divergent
formal power series partly because there are factorially many
Feynman diagrams, as was explained for example in \cite{GL1}. It is a
fundamental question how to give rigorous analytic meaning to these factorially
divergent series and how to numerically evaluate them for example for QCD,
and compare with well-tested experimental data. With this in mind,
't Hooft observed that when the size $N$ of the gauge is arbitrary,
the Feynman diagrams of perturbative gauge theory organize
themselves in ribbon graphs, i.e., abstract connected oriented surfaces of
some genus $g$ with some nonzero number of boundary components; see \cite{tH}.
Correspondingly, the free energy~$F_M(\tau,\hb)$ becomes a sum of power series
$F_{M,g}(\tau)$ of a single variable that carries the contribution of the
connected graphs of genus~$g$. 't Hooft suggested that

\begin{conjecture}
\label{conj.th}
There exists a disk $D_M$ that contains zero such that for all $g$,
$F_{M,g}(\tau)$ is a~power series analytic in $D_M$.
\end{conjecture}
If the above conjecture holds, then $F_{M,g}(\tau)$ can be evaluated
near $\tau=0$ by a convergent series. The above conjecture can be formulated
for any gauge theory. 't Hooft's conjecture has been verif\/ied in some simple
quantum-f\/ield theoretical models with~$U(N)$ symmetry, such as matrix models
(QFT in zero dimensions) and matrix quantum mechanics
(QFT in one dimension); see~\cite{BIPZ,BIZ,DGJ,EZ-Z,Wi}. Evidence for its
validity in $N=4$ super Yang--Mills theory has also appeared recently; see
\cite{BES,BSKS}.

\subsection{A free energy of a knot}
\label{sub.knot}

In the world of Quantum Topology, the free energy of a 3-manifold
is structurally similar to
two other invariants of knots: namely the {\em HOMFLY polynomial} $H_K(t,q)$
of a knot $K$ and the {\em Colored Jones function} $J_K(q,n)$, colored by the
$n$-dimensional irreducible representation of $\SL(2,\BC)$; see for
instance \cite{Tu}. The HOMFLY polynomial gives rise to an element of
$\hb^{-2} \BQ[[\tau,\hb^{2}]]$ using the substitution
\begin{gather}
\label{eq.tq}
t=e^{\tau}=e^{N \hbar}, \qquad q=e^{\hbar}.
\end{gather}
The same substitution works for the colored Jones function of a knot and is
the content of the so-called {\em loop expansion} of the colored Jones
function
\begin{gather}
\label{eq.loop}
J_K(q,n)=\sum_{n=0}^\infty F_{K,n}(\tau) \hbar^n \in \BQ[[\tau,\hbar]].
\end{gather}
The above loop expansion is really perturbation theory of Chern--Simons
theory along an Abelian f\/lat connection on the knot complement, see
\cite{Ro}.
For a detailed discussion on the meaning of the above expansion, see
\cite{Ro,GR,Ga1,GL2}. Moreover, its planar limit $F_{K,0}(\tau)$ can be
identif\/ied with the inverse Alexander polynomial $\D_K$
\begin{gather}
\label{eq.MMR}
F_{K,0}(\tau)=\frac{1}{\D_K(e^{\tau})}.
\end{gather}
This is the content of the Melvin--Morton--Rozansky Conjecture, shown in
\cite{B-NG}. In addition, for every $n \in \BN$, we have
\[
F_{K,n}(\tau)=\frac{P_{K,n}(e^\tau)}{\D_K(e^{\tau})^{2n+1}},
\]
where $P_{K,n}(t) \in \BZ[t^{\pm 1}]$ \cite{Ro}.
Since $\D_K(1)=1$, it follows that the radius of
convergence of~$F_{K,n}(\tau)$ at~$\tau=0$ is positive and independent of~$n$.

\subsection{Power series uniformly Gevrey-1}
\label{sub.gevrey}

In order to state our results, we need to formalize the analytic properties
of the free energy of a 3-manifold.
Recall that a formal power series
\[
f(\hb)=\sum_{n=0}^\infty a_n \hb^n
\]
is {\em Gevrey}-$s$ if there exists a positive constant $C>0$ such that
\[
|a_n| \leq C^n n!^s
\]
for all $n >0$.

\begin{definition}
\label{def.gevrey}
Consider a formal power series
\begin{gather}
\label{eq.fxe}
f(x,\e)=\sum_{n=-1}^\infty S_n(x) \e^n \in \e^{-1}\BC[[x,\e]]
\end{gather}
of two variables $(x,\e)$. We say that $f(x,\e)$ is
{\em Gevrey-$1$ with respect to $\e$, uniformly with respect to $x=0$}
(in short, $(x,\e)$ Gevrey-$1$) if there exists a constant $C>0$ so that
\begin{gather}
|[x^k]S_n(x)| \leq C^{n+k} n!
\end{gather}
for all $n,k \in \BN$. Here, $[x^k]g(x)$ denotes the coef\/f\/icient
of $x^k$ in a power series $g(x)$.
\end{definition}
Examples of power series $f(x,\e)$ $(x,\e)$ Gevrey-$1$ are the WKB solutions
of dif\/ference or dif\/ferential equations with a small parameter; see for
example \cite{AKKT}, where the authors call such series {\em uniformly
pre-Borel summable}. The loop expansion of a knot is $(\tau,\hbar)$
Gevrey-$1$; see \cite{Ga1}.

Observe that if a power series $f(x,\e)$ given in \eqref{eq.fxe} is
$(x,\e)$ Gevrey-$1$, then for all $n \in \BN$, the formal power series
$S_n(x)$ is analytic in a common neighborhood of $x=0$, independent of $n$.

\subsection{Statement of our results}
\label{sub.results}

For a closed 3-manifold $M$, let $Z_M$ denote the LMO invariant of $M$,
\cite{LMO}, which is an inf\/inite power series of vertex-oriented
trivalent graphs, and represents perturbation theory along the trivial
f\/lat connection of a 3-manifold. Given a metrized Lie algebra $\fg$, one
can replace graphs by rational numbers, keeping track of their number of
vertices, and thus create a power series $(W_{\fg} \circ Z_M)(\hb)
\in \BQ[[\hb]]$. Let $Z_M(N,\hbar)$ denote the LMO invariant
of $M$, composed with the $\gl_N$ weight system; see \cite{LMO}.
It is easy to see that $Z_M(N,\hbar)$ is a formal power series in $N$ and
$\hbar$ with constant term $1$, and its logarithm
\begin{gather}
\label{eq.FM2}
F_M(\tau,\hb)=\log Z_M(N,\hbar)
\end{gather}
which is by def\/inition the {\em free energy} of $M$,
can be written in the form~\eqref{eq.FM}, where $\tau=N\hbar$.

\begin{theorem}
\label{thm.1}
For every closed $3$-manifold $M$, the free energy $F_M(\tau,\hb)$
is $(\tau,\hb^{2})$ Gevrey-$1$.
\end{theorem}

Theorem~\ref{thm.1} follows from combining a presentation for the LMO
invariant given in Theorem~\ref{thm.2} with the def\/inition of the $\gl_N$
weight system, together with a crucial estimate on the number of
rooted maps in arbitrary genus obtained by Bender--Gao--Richmond following
work of Goulden--Jackson; see \cite{BGR,GJ} and Corollary \ref{cor.B} below.

\begin{theorem}
\label{thm.2}
For every closed $3$-manifold $M$, we can write its LMO invariant $Z_M$
in the form
\begin{gather}
\label{eq.Zma}
Z_M=\sum_{\Ga} c_{\Ga} \cdot \Ga,
\end{gather}
where we are summing over the set of trivalent graphs $\Ga$
and
\begin{gather}
\label{eq.cG}
|c_{\Ga}| \leq C_M^{n}
\end{gather}
for all $\Ga$ of degree $n$, where $C_M>0$ is a constant that depends on $M$.
\end{theorem}

Theorem \ref{thm.1} has the following corollary proving 't Hooft's conjecture.

\begin{corollary}
\label{cor.1}
For every closed $3$-manifold $M$, the power series $F_{M,g}(\tau)$
are analytic in a~common neighborhood of $\tau=0$, independent of $g$.
\end{corollary}

We now explain what happens when we specialize $N$ to be a f\/ixed natural
number, e.g.\ $N=2$. In \cite{GL1}, the f\/irst two authors proved the following
theorem for a f\/ixed metrized Lie algebra $\fg$. Let $Z_M$ denote the LMO
invariant of $M$ and let $W_{\mathfrak{g}}$ denote the corresponding weight
system. Then, $(W_{\mathfrak{g}} \circ Z_M)(\hbar) \in \BQ[[\hbar]]$.

\begin{theorem}[\cite{GL1}, Theorem~3]
\label{thm.GL}
For every metrized Lie algebra and every rational homology sphere $M$,
$(W_{\mathfrak{g}} \circ Z_M)(\hbar)$ is Gevrey-$1$.
\end{theorem}

\begin{corollary}
\label{cor.2}
Theorem {\rm \ref{thm.1}} implies Theorem {\rm \ref{thm.GL}} for $\mathfrak{g}=
\gl_N$ for every  fixed  $N$.
\end{corollary}

\subsection{Some calculations}
\label{sub.calc}

As a concrete illustration of our results, we can compute the free energy
of a Lens space. Let $L(d,b)$ denote the {\em Lens space} obtained by
$d/b \in \BQ$ surgery on the unknot in $S^3$.

Recall the $\a$-polylogarithm function
\begin{gather}
\label{eq.polylog}
\Li_{\a}(x)=\sum_{n=1}^\infty \frac{x^n}{n^{\a}}
\end{gather}
for $\a \in \BR$, def\/ined by the absolutely convergent series for $|x|<1$
and analytically continued in $\BC\setminus\{0,1\}$. For a detailed discussion,
see~\cite{Oe} and also~\cite{CG}. Let $B_n$ denote the $n$th
{\em Bernoulli number}
def\/ined by the generating series
\[
\frac{x}{e^x-1}=\sum_{n=0}^\infty \frac{B_n}{n!} x^n.
\]

\begin{theorem}
\label{thm.lens}
Consider a Lens space $M=L(d,b)$. Its free energy is given by
\begin{eqnarray}
\label{eq.Flens2}
F_{M,g}(\tau) &=&
(2g-1)\frac{B_{2g}}{(2g)!} \left(d^{2-2g}\Li_{3-2g}(e^{\tau/d})-
\Li_{3-2g}(e^{\tau}) \right) +a_g(\tau),
\end{eqnarray}
where
\begin{gather}
\label{eq.ag}
a_g(\tau)=\begin{cases}
\displaystyle -\frac{\tau^2}{2} \log d - (d^2-1) \zeta(3)+ \l_{L(d,b)} \frac{\tau^3}{2}
& \text{if} \ \  g=0, \vspace{1mm}\\
\displaystyle \frac{\tau}{24}(1-d^{-1}) + \frac{1}{12}\log d - \l_{L(d,b)} \frac{\tau}{2}
& \text{if} \ \  g=1, \vspace{1mm}\\
0 & \text{if} \ \  g \geq 2.
\end{cases}
\end{gather}
\end{theorem}
It follows that for all $g \geq 2$ (resp.~$g=0,1$),
$F_{M,g}$ has analytic continuation as a meromorphic (resp.\ multivalued
analytic) function on $\BC\setminus \frac{1}{d\BZ^*(1)}$
(resp. $\BC\setminus \frac{1}{d\BZ(1)}$), where
\begin{gather}
\label{eq.BZ1}
\BZ(1)=2 \pi i \BZ, \qquad \BZ^*(1)=\BZ(1)\setminus\{0\}.
\end{gather}

\section[The LMO invariant and its $\ell^\infty$ Gromov norm]{The LMO invariant and its $\boldsymbol{\ell^\infty}$ Gromov norm}
\label{sec.LMO}

\subsection[The $\ell^\infty$ Gromov norm]{The $\boldsymbol{\ell^\infty}$ Gromov norm}
\label{sub.LMO}

The purpose of this section is to prove Theorem~\ref{thm.2} extending
the Gromov norm techniques from~\cite{GL1}. We will assume some familiarity
with~\cite{GL1}. To begin with, the LMO invariant~$Z_M$ of a closed 3-manifold
takes value in a completed graded vector space~$\calA(\emptyset)$ of
vertex-oriented trivalent graphs, modulo some linear homogeneous
$\mathrm{AS}$ and $\mathrm{IHX}$
relations, where the degree of a graph is half the number of vertices.
The LMO invariant gives a meaningful def\/inition to the Chern--Simons
perturbation
theory along a trivial f\/lat connection, and the trivalent graphs mentioned
above as the Feynmann diagrams of a~$\phi^3$ gauge theory with Chern--Simons
action. For a detailed discussion, see \cite{KT}.
 The linear $\mathrm{AS}$ and $\mathrm{IHX}$ relations are the
diagrammatic version of the antisymmetry and the Jacobi identity of the
Lie algebra of the gauge theory.

Let $\calA_n(\emptyset)$ denote the subspace of $\calA(\emptyset)$ of degree
$n$. Then, $\calA_n(\emptyset)$ is a f\/inite dimensional vector space spanned
by the f\/inite set of vertex-oriented trivalent graphs with $2n$ vertices.
Of course,
this spanning set is not a basis, due to the linear $\mathrm{AS}$ and
$\mathrm{IHX}$ relations. This motivates the following concept of the
$\ell^p$ Gromov norm, extending the one in \cite[Def\/inition~1.1]{GL1}.

\begin{definition}
\label{def.gnorm}
Consider a vector space $V$ and a spanning set $b$. Fix $p \in [1,\infty]$.
For $v \in V$, def\/ine the $\ell^p$ {\em norm} by:
\begin{gather}
\label{eq.gnorm}
|v|_p=\begin{cases}
\inf \big(\sum_j |c_j|^p\big)^{1/p} & \text{if} \ \ p \in [1,\infty), \\
\inf \max_j |c_j|           & \text{if} \ \ p=\infty,
\end{cases}
\end{gather}
where the inf\/imum is taken over all presentations of the form
$v=\sum_j c_j v_j$, $v_j \in b$.
\end{definition}

We now apply the above def\/inition to each of $\calA_n(\emptyset)$ and combine
them into a single power series.

\begin{definition} \label{def.gnormA}\quad {}

\begin{enumerate}\itemsep=0pt
\item[\rm{(a)}] Fix $p \in [1,\infty]$ and $n \in \BN$ and consider the vector space
$\calA_n(\emptyset)$ spanned by the set $b$ of oriented trivalent graphs of
degree $n$. For $v \in \calA_n(\emptyset)$, we denote by $|v|_p$ the norm of~$v$.

\item[\rm{(b)}] If $v \in \calA(\emptyset)$, we def\/ine
\begin{gather}
\label{eq.gnormA}
|v|_p(\hb) =\sum_{n=0}^\infty |\pi_n(v)|_p \hb^n \in \BQ[[\hb]]
\end{gather}
where $\pi_n: \calA(\emptyset) \longto \calA_n(\emptyset)$ denotes the
projection in $\calA_n(\emptyset)$.

\item[\rm{(c)}] We say that $v \in \calA(\emptyset)$ has Gevrey-$s$ $p$-norm if
$|v|_p(\hb)$ is Gevrey-$s$.
\end{enumerate}
\end{definition}

\subsection[A brief review of \cite{GL1}]{A brief review of \cite{GL1}}
\label{sub.reviewGL}

The main result of \cite{GL1} was the following theorem.

\begin{theorem}[\cite{GL1}, Theorem~2]\label{thm.g1}
For every integral homology sphere $M$, the power series $|Z_M|_1(\hb)
\in\BQ[[\hb]]$ is Gevrey-$1$.
\end{theorem}
Theorem~\ref{thm.g1} is a convergence property of a {\em universal finite type
invariant} of homology spheres, and uses a well-known universal
f\/inite type invariant of links in $S^3$, namely the Kontsevich integral.
For an introduction to the notion of f\/inite type invariants of knots and
3-manifolds, see~\cite{B-N,Oh2}. Since our proof of Theorem~\ref{thm.2} will
use some of the ideas of the proof of Theorem~\ref{thm.g1}, we begin by
recalling the proof of Theorem~\ref{thm.g1}, which goes as follows.
First, we begin with
the following description of the LMO invariant from~\cite{BGRT}.

\begin{itemize}\itemsep=0pt
\item
Fix a presentation of an integral homology sphere $M$ as surgery on a
unit-framed boundary link $L$ in $S^3$.
\item
Choose a presentation of $L$ as the closure of a framed string link $T$.
\item
Consider the normalized Kontsevich integral $\check Z_T$, which takes values
in the completed $\BQ$-vector space $\calA(\star_X)$ of vertex-oriented
Jacobi diagrams with legs colored by a set $X$ in 1-1 correspondence with
the components of $T$.
\item
Perform {\em formal diagrammatic Gaussian integration}
\[
\int  dX: \calA(\star_X) \longto \calA(\emptyset)
\]
to $\check Z_T$. I.e., separate out the {\em strut part} of $\check Z_T$ (which
is an exponential of Jacobi diagrams with no trivalent vertices)
and join the
legs of the remaining diagrams using the inverse linking matrix of the
strut part.
\item
Finally, def\/ine
\begin{gather}
\label{eq.LMOdef}
Z_M=\frac{\int \check Z_L dX}{(\int \check Z_{U^+} dX)^{\s_+(L)}
(\int \check Z_{U^-} dX)^{\s_-(L)}},
\end{gather}
where $U^{\pm}$ denotes the $\pm 1$ framed unknot in $S^3$ and
$\s_{\pm}(L)$ denotes the number of positive (resp. negative)
eigenvalues of the linking matrix of $L$.
\end{itemize}

Next, we extend our notion of $\ell^p$ norm to all intermediate spaces of
Jacobi diagrams (with or without skeleton, and with or without symmetrized
legs).

Next, we show that the Kontsevich integral is Gevrey-$0$.

\begin{theorem}[\cite{GL1}, Theorem~8]
\label{thm.Ko}
For every framed tangle $T$, $|Z_T|_1(\hb) \in \BR[[\hb]]$ is Gev\-rey-$0$.
\end{theorem}

This follows from the def\/inition of the Kontsevich integral using the
KZ associator, together with the following lemma.

\begin{lemma}[\cite{GL1}, Proposition~2.6]
\label{lem.unknot}
If $U$ denotes the zero-framed unknot, then $|Z_U|_1(\hb) \in \BR[[\hb]]$
is Gevrey-$0$.
\end{lemma}

Finally we use the following key lemma.

\begin{lemma}[\cite{GL1}, Lemma~2.10]
\label{lem.gl}
For every boundary string link $T$ with framing $\pm 1$,
\linebreak $| \int  \check Z_T dX|_1(\hb)$ $\in \, \BR[[\hb]]$ is Gevrey-$1$.
\end{lemma}

The condition on the string link being a boundary one is to ensure that
the strutless part of~$\check Z_T$ contains no trees; see~\cite{HM}. Lemma~\ref{lem.gl} takes care of the normalization of the LMO invariant given
by~\eqref{eq.LMOdef}. Ignoring technicalities,
the main point in the proof of Lemma~\ref{lem.gl} is the following.

If $G$ is uni-trivalent graph of degree $n+k$ (i.e., with $2n+2k$ vertices)
with $2k$ legs, then there are $(2k-1)!!=1\cdot 3\cdot 5\cdot \cdots \cdot (2k-1)$ ways to pair
the legs of $G$ together. If $G$ has no tree components, then $k \leq n$
thus
\[
(2k-1)!! \leq (2n-1)!! \leq C^n n!.
\]

\begin{remark}
\label{rem.1p}
If we replace $|\cdot|_1$ by $|\cdot|_p$ for $p \in [1,\infty)$, Theorems~\ref{thm.g1},~\ref{thm.Ko} and Lemmas~\ref{lem.unknot},~\ref{lem.gl}
remain true.
\end{remark}

In the next section we discuss what happens when we use the $|\cdot|_\infty$
rather than the $|\cdot|_1$ norm.

\subsection{Proof of Theorem \ref{thm.2}}
\label{sub.thm2}

We now show how to modify the statements of the previous section in order
to show the following reformulation of Theorem \ref{thm.2}. In that theorem
we need to f\/ix a closed 3-manifold. When $M$ has positive Betti number, the LMO
invariant simplif\/ies a lot, and it is possible to give a direct proof by
an explicit formula. For example, when the f\/irst Betti number $b_1(M)$
is greater than~$3$, then the LMO invariant vanishes and Theorem \ref{thm.2}
obvious holds. When $1< b_1(M) \leq 3$, the LMO invariant can be directly
computed, see~\cite{HT,Oh2}. When $b_1(M)=1$, the LMO invariant is computed by~\cite{GH} in terms of the Alexander polynomial of~$M$, and Theorem~\ref{thm.2} follows immediately. Thus, it suf\/f\/ices to assume that $M$ is
a rational homology sphere, and further (using the multiplicative behavior
of the LMO with respect to connected sums), it suf\/f\/ices to assume that $M$ is
obtained by surgery on a rationally framed boundary link in~$S^3$. In that
case, the LMO invariant can be computed via the Aarhus integral.

\begin{theorem}
\label{thm.2alt}
For every rational homology sphere $M$, $|Z_M|_{\infty}(\hb) \in \BR[[\hb]]$
is Gevrey-$0$.
\end{theorem}

This theorem follows from the following.

\begin{theorem}
\label{thm.Koalt}
For every framed tangle $T$, $|Z_T|_\infty(\hb) \in \BR[[\hb]]$ is Gevrey-$0$.
\end{theorem}

This follows easily from Theorem \ref{thm.Ko} and the fact that $|v|_\infty
\leq |v|_1$ for all $v \in \calA_n(\emptyset)$.

\begin{lemma}
\label{lem.glalt}
For every framed string link $T$ with invertible linking matrix
$|\int \check Z_T dX|_\infty(\hb) \in \BR[[\hb]]$ is Gevrey-$0$.
\end{lemma}

The proof of this lemma is as follows. First we write
$\check Z_T$ as follows, using Theorem \ref{thm.Koalt}
\[
\check Z_T=\exp\left( \frac{1}{2} \sum_{ij} l_{ij} |_i^j\right) \sum_{\Ga} c_{\Ga}
\cdot \Ga,
\]
where $(l_{ij})$ is the linking matrix of $T$, $|_i^j$ is a strut colored
by the components $i$ and $j$ of $T$, the summation is over the set
of uni-trivalent graphs with no strut components and legs colored by the
components of $T$ and
\[
|c_{\Ga}| \leq C^{\deg(\Ga)}
\]
for some constant $C>0$. If $G$ is a trivalent graph of degree
$n$, it has $3n$ edges. If $G$ comes from pairing the legs of a uni-trivalent
graph $\Ga$, then $\Ga$ is obtained from $G$ by cutting some (let's say $k
\leq 3n $) of the edges of $G$ in half. Such a uni-trivalent graph $\Ga$ has
degree $n+k$ (i.e., $2n+2k$ vertices) and $2k$ legs and its coef\/f\/icient in
$\check Z_T$ is bounded by $C^{n+k} \leq C^{4n}$.
Since~$G$ has~$3n$ edges, there are at most $2^{3n}$ such graphs $\Ga$.
Thus, pairing the legs of the strutless part of $\check Z_T$, we obtain
that
\[
\int \check Z_T dX =  \sum_{G} c_{G}
\cdot G,
\]
where
\[
|c_G| \leq C'^{\deg(G)}
\]
for some $C'>0$. This concludes the proof of Theorem~\ref{thm.2}.

\section[The $\gl_N$ weight system]{The $\boldsymbol{\gl_N}$ weight system}
\label{sec.glN}

The LMO invariant of a closed 3-manifold takes values in a completed
vector space $\calA(\emptyset)$ spanned by trivalent graphs. To get a
numerical invariant with values in $\BQ[[\hb]]$, we need to replace every
graph by a combinatorial weight. This is exactly, what a weight system does.
More precisely, given a Lie algebra $\fg$ with an invariant inner product,
there is a weight system $\BQ$-linear map:
\[
W_{\fg}: \calA(\emptyset) \longto \BQ[[\hb]]
\]
discussed at great length in~\cite{B-N}. In \cite{B-N}, the following
description of $W_{\fg}$ is given.
\[
W_{\gl_N}: \calA(\emptyset) \longto \BQ[[N,\hb]]
\]
is def\/ined by:
\begin{gather}
\label{eq.wglN}
D \longto \sum_M (-1)^{s_M} N^{b_{D,M}} \hb^{2g_{D,M}-2+b_{D,M}},
\end{gather}
where
\begin{itemize}\itemsep=0pt
\item
the sum is over all possible markings $M$ of the trivalent vertices
of $D$ by $0$ or $1$,
\item
$s_M$ is the sum, over the set of trivalent vertices, of the values of $M$,
\item
$\S_{D,M}$ denotes the $X$-marked surface
obtained by thickening the trivalent vertices of $D$ as follows:

\begin{figure}[!ht]
\centerline{\includegraphics[scale=0.6]{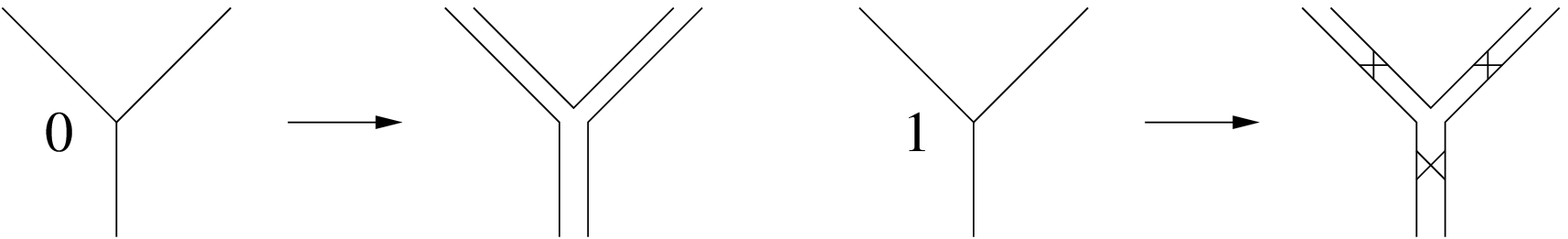}}
\label{figure}
\end{figure}

and thickening the edges of $D$, and connecting up into to a surface.
It turns out that $\S_{D,M}$ is well-def\/ined, oriented and perhaps
disconnected,
\item
$g_{D,M}$ and $b_{D,M}$ denote the genus (i.e., sum of the genera of
its connected components) and the number of boundary components
of $\S(D,M)$,
\item
it turns out that the degree of $D$ is $2g_{D,M}-2+b_{D,M}$.
\end{itemize}

For example, we have
\[
W_{\gl_N}(\Theta)=2(N^3-N)\hb,
\]
where $\Theta$ is the obvious planar trivalent graph with counterclockwise
cyclic order at each vertex. The next lemma follows directly from the above
def\/inition of the $W_{\gl_N}$ weight system.

\begin{lemma}
\label{lem.WglN}
If $\Ga$ is a connected trivalent graph of degree $n$, then
\begin{gather}
\label{eq.WglN}
W_{\gl_N}(\Ga)= p_{\Ga}(N) \hbar^n,
\end{gather}
where $p_{\Ga}(N) \in \BZ[N]$ is a polynomial in $N$ of degree at most $n+2$
and $\ell^1$-norm at most $2^n$.
\end{lemma}

\section{Asymptotics of the number of rooted maps of arbitrary genus}
\label{sec.comb}

In a series of papers in various collaborations, Bender et al studied the exact
and asymptotic number of rooted maps on a surface; see \cite{BC,BGR}.
Recall that a {\em map} $(G,S)$ is a graph $G$ embedded in a connected,
oriented, closed topological surface $S$ in such a way that each maximal
connected component of $S\setminus G$ is a topological disk. A map is
{\em rooted} if an edge, a direction along the edge, and a side of the edge
are distinguished. Let $T_g(n)$ denote the number of $n$-edged rooted maps
on a surface of genus $g$. In \cite{BC}, Bender--Canf\/ield, using generating
functions, and the Darboux--Polya enumeration method,  give the following
asymptotic expansion of $T_g(n)$.

\begin{theorem}[\cite{BC}]
\label{thm.B1}
For $g$ fixed and $n \longto \infty$ we have
\begin{gather}
\label{eq.Tgn}
T_g(n) \sim t_g n^{5(g-1)/2} 12^n,
\end{gather}
where the constants $t_g$ are computable via computable non-linear recursions.
In particular,
\[
t_0=\frac{2}{\sqrt{\pi}}, \qquad t_1=\frac{1}{24}, \qquad
t_2=\frac{7}{4320 \sqrt{\pi}}.
\]
\end{theorem}

Many interesting families of maps satisfy asymptotic formulas
of the form
\[
\a t_g (\b n)^{5(g-1)/2} \ga^n.
\]
This is discussed extensively for example in \cite[Section~4]{Wi}; see also
\cite{G1,G2}.
For over 20 years, the constants $t_g$ remained hard to compute or estimate,
partly due to the complexity of their non-linear recursion.
A breakthrough was achieved recently. Using work of Goulden--Jackson~\cite{GJ}
on the KP hierarchy, Bender--Gao--Richmond simplif\/ied the non-linear
recursion relation for~$t_g$, and obtained the following.

\begin{theorem}[\cite{BGR}]
\label{thm.B2}
We have
\begin{gather}
\label{eq.tg}
t_g \sim \frac{40 \sin(\pi/5) K}{\sqrt{2 \pi}}
\left(\frac{1440g}{e}\right)^{-g/2},
\end{gather}
where $K=0.10486898772254091800\dots$ is a constant.
\end{theorem}
The numerical value of $K$ was given incorrectly at \cite{BGR}. An exact
value of $K$ is given in equation~\eqref{eq.Kexact} in
Section~\ref{sec.appendix} below.
For our purposes, it suf\/f\/ices to note that a trivalent graph of degree $n$
has $3n$ edges. Taking into account the choice of a root, it follows that

\begin{corollary}
\label{cor.B}
There are at most $12n T_g(3n)$ connected trivalent graphs of genus $g$,
where $T_g(n)$ satisfies~\eqref{eq.Tgn} and~\eqref{eq.tg}.
\end{corollary}

\begin{remark}
\label{rem.unoriented}
Theorem \ref{thm.B1} has an analogous statement (with a dif\/ferent constant
$p_g$ instead of $t_g$) for rooted graphs in unoriented surfaces; see
\cite[Theorem~1]{BC}. Without doubt, there is an asymptotic expansion for
$p_g$ analogous to $t_g$. In fact, one can write down explicitly a non-linear
recursion relation that $p_g$ ought to satisfy, and study its exact asymptotic
behavior. This will be discussed in a separate publication; see~\cite{GM3}.
\end{remark}

\section{Proofs}
\label{sec.proofs}

\subsection{Proof of Theorem~\ref{thm.1}}
\label{sec.thm1}

Let $Z_M(N,\hbar)$ denote the LMO invariant of $M$, composed with the weight
system of $\gl_N$. Then, we can write
\begin{gather}
\label{eq.ZM}
\log Z_M(N,\hbar)=\sum_{2g-2+d>0, d>0} a_{M,g,d} N^d \hbar^{2g-2+d}.
\end{gather}
Recall that the free energy of $M$ is def\/ined by \eqref{eq.FM2}
\[
F_M(\tau,\hb)=\log Z_M(N,\hbar),
\]
where $\tau=N\hbar$. It follows that the free energy has the form
\eqref{eq.FM} where
\begin{gather}
\label{eq.FMg}
F_{M,g}(\tau)=\sum_{d: 2g-2+d>0, d>0} a_{M,g,d} \tau^d \in \BQ[[\tau]].
\end{gather}

On the other hand, Theorem \ref{thm.2} implies that we can write
\[
\log Z_M=\sum_{\Ga} c_{\Ga} \cdot \Ga,
\]
where the sum is over a set of connected trivalent graphs, where there
exists a constant $C$ such that
\begin{gather}
\label{eq.cGa}
|c_{\Ga}| \leq C^n
\end{gather}
if $\Ga$ has degree $n$. Applying the $W_{\gl_N}$ weight system,
and using Lemma \ref{lem.WglN} it follows that
\[
F_M(\tau,\hb)= \sum_{n=1}^\infty \sum_{\Ga} c_{\Ga} \cdot p_{\Ga}(N) \hbar^n,
\]
where the $\Ga$ summation is over a set of connected trivalent graphs
of degree $n$. Thus,
\begin{gather}
\label{eq.aMgd}
a_{M,g,d}=\sum_{\Ga} c_{\Ga} \cdot [N^d]p_{\Ga}(N),
\end{gather}
where the sum is over the set of connected trivalent graphs $\Ga$ of degree
$n=2g-2+d$ that embed on a surface of genus $g$. Lemma~\ref{lem.WglN}
estimates the coef\/f\/icients $[N^d]p_{\Ga}(N)$ and since $2g-2+d=n$,
it implies that
\[
|a_{M,g,d}| \leq  C_1^n \sum_{\Ga} |c_{\Ga}|
\]
for some constant $C_1>0$.
Equation \eqref{eq.cGa} together with $2g-2+d=n$ implies that
\[
|a_{M,g,d}| \leq  C_2^n \sum_{\Ga} 1
\]
for some constant $C_2>0$. The above sum is over the set of trivalent graphs
of degree $n=2g-2+d$ (and thus, $3n$ edges) that embed on a surface of genus
$g$. Corollary \ref{cor.B} implies that
\begin{gather*}
\sum_{\Ga} 1  \leq   8n T_g(3n)   \leq  C_1^n \left(\frac{(3n)^5}{g}\right)^{g/2}
 \leq  C_2^n \frac{(2g-2+d)^{5g/2}}{g^{g/2}} \\
 \phantom{\sum_{\Ga} 1}{} =  C_2^n \frac{(2g-2+d)^{5g/2}}{(2g)^{5g/2}} \frac{(2g)^{5g/2}}{g^{g/2}}.
\end{gather*}
Since $(1+a/n)^n \leq e^a$, we have
\[
\frac{(2g-2+d)^{5g/2}}{(2g)^{5g/2}} \leq e^{5(d-2)/4}.
\]
Moreover,
\[
\frac{(2g)^{5g/2}}{g^{g/2}} \leq C_3^g (2g)!
\]
Since $n=2g-2+d$, it follows that
\[
|a_{M,g,d}| \leq C_3^{g+d} (2g)!
\]
and concludes the proof of Theorem~\ref{thm.1}.

\subsection{Proof of Corollary \ref{cor.2}}
\label{sub.cor2}

Fix a natural number $N \in \BN$. It suf\/f\/ices to prove that
$(W_{\gl_N} \circ \log Z_M)(\hbar) \in \BQ[[\hbar]]$ is Gevrey-$1$.
Equation \eqref{eq.ZM} implies that
\begin{gather}
[\hbar^n] (W_{\gl_N}\circ\log Z_M)(\hbar) =
\sum_{g,d: 2g-2+d=n} a_{M,g,d} N^d.
\end{gather}
The proof of Theorem \ref{thm.1} implies that there exists $C>0$ so that
\[
|a_{M,g,d}| \leq C^{g+d} (2g)!.
\]
Since $2g-2+d=n$, $C^{g+d} \leq C_1^{3n/2}=C_2^n$ and $(2g)! \leq C_3^n n!$.
Thus,
\begin{gather*}
|\coeff(\hbar^n, (W_{\gl_N}\circ\log Z_M)(\hbar))|
 \leq  C_4^n n! \sum_{g,d: 2g-2+d=n} N^d  \leq  C_4^n n! N^{n+2} (n+2).
\end{gather*}
Since $N$ is {\em fixed}, it follows that
\[
|[\hbar^n] (W_{\gl_N}\circ \log Z_M)(\hbar) | \leq C^n n!
\]
which concludes the proof of Corollary~\ref{cor.2}.

\section{The free energy of a Lens space}
\label{sec.calc}

In this section we will prove Theorem \ref{thm.lens}.
The next proposition computes the image of the
LMO invariant of a Lens space $L(d,b)$ under the weight system $W_{\fg}$.

Let $\fg$ denote a metrized Lie algebra with inner product $\la \cdot,\cdot\ra$
and let $\Phi_+$ denote the positive roots of $\fg$, $|\Phi_+|$ denote
the cardinality of $\Phi_+$  and $\rho$ denote
half the sum of the positive roots of~$\fg$. Let $c_{\fg}$ denote the
product of the quadratic Casimir of $\fg$ with the dimension of $\fg$.
Let $\lambda_M$ denote Casson invariant.

\begin{proposition}
\label{prop.Lens}
With the above assumptions we have
\begin{gather}
\label{eq.Lens}
(W_{\fg} \circ Z_{L(d,b)})(\hb)=
\exp\left(\frac{\l_{L(d,b)}}{4}c_{\fg} \hb\right) d^{|\Phi_+|}
\prod_{\a \in \Phi_+}
\frac{\sinh\left((\a,\rho)\hb/(2d)\right)}{
\sinh\left((\a,\rho)\hb/2\right)}
\in \BQ[[\hb]].
\end{gather}
\end{proposition}

\begin{proof}
There are two proofs of this proposition. One proof uses
\begin{itemize}\itemsep=0pt
\item[(a)] The existence of the Ohtsuki series (which come from the
Reshetikhin--Tuarev invariants of 3-manifolds). This was shown
by Ohtsuki in~\cite{Oh1} for $\PSU(2)$, and by the second author in~\cite{L1}
for $\PSU(N)$ and more generally in \cite{L2} for all simple Lie algebras.
\item[(b)] The computation of the Ohtsuki series for Lens spaces, done
by Takata in~\cite{Ta}.
\item[(c)]  The identif\/ication of the left hand side of~\eqref{eq.Lens} with
the Ohtsuki series. This was the result of the unpublished fourth part of~\cite{BGRT} and independently an unpublished work of the second author that
was recently completed by Kuriya~\cite{Ku}, see also~\cite{Oh2}.
\end{itemize}
An alternative proof uses the computation of the LMO invariant for Lens spaces
by Bar-Natan and Lawrence \cite{B-NL}. Let us give the details of this proof.
We will use the notation from \cite{B-NL}.
In \cite[Proposition~5.1]{B-NL} it is shown that the LMO invariant of $L(d,b)$
is given by:
\[
Z_{L(d,b)}=\exp\left(\frac{-S(b/d)}{48}\Theta\right)
\frac{\la \Om_x, \Om_{x/d} \ra_x}{\la \Om_x, \Om_{x} \ra_x}.
\]
Moreover, $-S(b/d)=12 \l_{L(b,d)}$.
Now, we apply the weight system $W_{\fg}$. We have
\begin{gather*}
W_{\fg}(\Theta)(\hb) = c_{\fg} \hb, \\
W_{\fg}(\la \Om_x, \Om_{x} \ra_x)(\hb) =
\prod_{\a \in \Phi_+}\frac{\sinh\left((\a,\rho)\hb/2\right)}{
(\a,\rho)\hb/2}
\end{gather*}
and
\[
W_{\fg}(\la \Om_x, \Om_{x/d} \ra_x)(\hb) =
W_{\fg}(\la \Om_x, \Om_{x} \ra_x)(\hb/d).
\]

The result follows.
\end{proof}

\subsection{Some special functions}
\label{sub.special}

Let us introduce some auxiliary functions which have already appeared in the
LMO invariant of the Lens spaces and which are important ingredients of the
proof of Theorem \ref{thm.lens}. Consider the functions $f$ and $\ti f$
def\/ined by
\[
\tilde f(x):= \frac{\sinh(x/2)}{x/2}, \qquad \text{and} \qquad
f(x) := \log \tilde f(x).
\]
It is known that
\[
f(x) = \sum_{k=1}^\infty b_{2k} x^{2k} = \frac{x^2}{24}-\frac{x^4}{2880}+
\frac{x^6}{181440} + \cdots,
\]
where $b_{2k}$'s are the modif\/ied Bernoulli numbers and are related to the
ordinary Bernoulli numbers~$B_k$ as follows:
\[
b_{2k} = \frac{B_{2k}}{2k (2k)!}.
\]
The f\/irst polylogarithm is given by
\begin{gather}
\label{eq.Li1}
\Li_1(x)=-\log(1-x).
\end{gather}
The relation between $f(x)$ and $\Li_1$ is the following
\begin{gather}
\label{eq.L1}
f(x)=-\Li_1(e^{x})-\frac{x}{2}-\log(-x).
\end{gather}
The left hand side is an analytic function in the unit disk $|x| < 1$ with
an analytic continuation as a multivalued analytic function in $\BC\setminus
2 \pi i (\BZ\setminus \{0\})$. The right hand side is well-def\/ined in the
unit disk minus $[0,1)$, and extends as a continuous function over the cut
$[0,1)$. Then, both sides agree on the unit disk $|x|<1$.
Equation \eqref{eq.L1} follows from the following easy computation
\[
f(x) = \log\left(\frac{e^{x/2}-e^{-x/2}}{x}\right) =
\log\left( \frac{e^{-x/2}}{-x} (1-e^x)\right)
\]
and equation \eqref{eq.Li1}.

From the power series def\/inition of the polylogarithm, it is easy to see that
for all $\a \in \BR$ we have
\begin{gather}
\label{eq.derLi}
\frac{d}{dx} \Li_{\a}(e^x)=\Li_{\a-1}(e^x).
\end{gather}
Moreover, when $\a$ is a negative integer, then
\begin{gather}
\label{eq.Lin}
\Li_{\a}(x)=\frac{P_{\a}(x)}{(1-x)^{-\a+1}},
\end{gather}
where $P_{\a}(x) \in \BZ[x]$ is a palindromic polynomial of degree $-\a$.
It follows that for $\a$ a negative integer, the function $\Li_{\a}(e^\tau)$
is an even meromorphic function on $\BC$ with poles at $\BZ(1)$, where the
latter set is def\/ined in~\eqref{eq.BZ1}.

\subsection{Proof of Theorem \ref{thm.lens}}
\label{sub.thm.lens}
This section is devoted to the proof of Theorem \ref{thm.lens}. Let us observe
that the $\gl_N$ and the $\sl_N$ weight system agree on all nonempty trivalent
graphs (since the structure constants of an Abelian Lie algebra vanish,
and are placed at the trivalent vertices of a trivalent graph),
and that the logarithm of the LMO invariant of a closed 3-manifold
is a series of nonempty connected trivalent graphs. Therefore, to compute
the free energy, we can work with the $\sl_N$ weight system.
Recall that the set $\Phi_+$ of positive roots of $\mathfrak{sl}_N$ has
$N(N-1)/2$ elements. From Proposition \ref{prop.Lens} one has
\[
(W_{\mathfrak{sl}_N} \circ Z_{L(d,b)})(\hb)=
 e^{\frac{\l_{L(d,b)}}{12} c_{\sl_N}\hb}
\prod_{\a \in \Phi_+}
\frac{\tilde f ((\a, \rho)\hb/d)}{\tilde f( (\a, \rho)\hb)}.
\]
Taking the logarithm, we get
\[
F_M(\tau,\hb) = \hb\, \frac{\l_{L(d,b)}}{12} c_{\mathfrak{sl}_N}
+ \sum_{\a \in \Phi_+}\left( f ((\a, \rho)\hb/d) - f ((\a, \rho)\hb)\right).
\]
Note that
\[
\hb\,
\frac{\l_{L(d,b)}}{4} c_{\mathfrak{sl}_N}=\frac{\l_{L(d,b)}}{2}\hb(N^3-N)=
\frac{\l_{L(d,b)}}{2}\left(\frac{\tau^3}{\hb^2}-\tau\right).
\]
When $\a$ runs through the set $\Phi_+$, $(\a,\rho)$ takes the value
$j$ the total of $N-j$ times, for every $j=1,\dots,N-1$. Hence
\begin{gather}
\label{eq.Ftemp3}
F_M(\tau,\hb) = \hb\, \frac{\l_{L(d,b)}}{12} c_{\fg}
+ \sum_{j=1}^{N-1} \left( (N-j) f (j\hb/d) - (N-j) f (j\hb)\right).
\end{gather}

Now we have
\begin{gather*}
\sum_{j=1}^{N-1} (N-j) f(j\hb)   = \sum_{j=1}^{N-1} (N-j)
\sum_{k=1}^\infty  b_{2k} j^{2k} \hb^{2k}
 = \sum_{k=1}^\infty  b_{2k} \hb^{2k}  \sum_{j=1}^{N-1} (N-j) \, j^{2k}
\end{gather*}

Using the well-known identity expressing the sum of powers by a
Bernoulli polynomial
\[
\sum_{j=1}^{N-1} j^k = \frac{1}{k+1} \sum_{s=0}^k \binom{k+1}{s}
B_s N^{k+1-s},
\]
it follows that
\[
\sum_{j=1}^{N-1} (N-j) \, j^{2k}=
\sum_{g=0}^k  \frac{(2k)! \, (1-2g)}{(2g)!(2k+2-2g)!} B_{2g} N^{2k+2-2g}.
\]

Therefore, using $N=\tau/\hb$, we get
\begin{gather}
\sum_{j=1}^{N-1} (N-j) f(j\hb)
 = \sum_{k=1}^\infty  b_{2k} \hb^{2k}
\left \{  \sum_{g=0}^k  \frac{(2k)! \,
(1-2g)}{(2g)!(2k+2-2g)!} B_{2g} N^{2k+2-2g} \right\}\notag \\
\phantom{\sum_{j=1}^{N-1} (N-j) f(j\hb)}{}
= \sum_{g=0}^\infty \frac{(1-2g)\, B_{2g} \hb^{2-2g}}{(2g)!}
\sum_{k=\max\{g,1\}}^\infty
\frac{(2k)!}{(2k+2-2g)!}b_{2k} \tau^{2k+2-2g}\notag \\
\phantom{\sum_{j=1}^{N-1} (N-j) f(j\hb)}{}
= \sum_{g=0}^\infty \frac{(1-2g)\, B_{2g} \hb^{2-2g}}{(2g)!}
\sum_{l=0}^\infty
\frac{(2l+2g)!}{(2l+2)!}b_{2g+2l} \tau^{2l+2}\label{eq.Ftemp2}
\end{gather}
with a minor variation when $g=0$. Let us def\/ine the auxiliary functions
$F_g$ by
\begin{gather}
\label{eq.Fglens}
F_g(\tau):= \sum_{l=0}^\infty \frac{(2l+2g)!}{(2l+2)!}b_{2g+2l} \tau^{2l+2}.
\end{gather}
Then, equations \eqref{eq.Ftemp3}, \eqref{eq.Ftemp2} and \eqref{eq.Fglens}
and the replacement of $(\hb,\tau)$ by $(\hb/d,\tau/d)$ imply that
\begin{gather}
\label{eq.Ftemp7}
F_{M,g}(\tau)=
(1-2g)\frac{B_{2g}}{(2g)!} \left(d^{2-2g} F_g(\tau/d)- F_g(\tau)
\right) + \frac{\l_{L(d,b)}}{2}(\tau^3 \d_{g,0}-\tau \d_{g,1}).
\end{gather}
We claim that the functions $F_g$ are given by
\begin{gather}
\label{eq.Ftildeg}
F_g(\tau)=-\Li_{3-2g}(e^\tau) +
\begin{cases}
\displaystyle -\frac{\tau^2}{2}\log(-\tau)-\frac{\tau^3}{12}+\frac{3\tau^2}{4}
-\frac{\pi^2\tau}{6}+\zeta(3)
  & \text{if} \ \  g=0 , \vspace{1mm}\\
\displaystyle  -\frac{\tau}{2}-\log(-\tau) & \text{if} \ \  g=1,\vspace{1mm} \\
\displaystyle (2g-3)! \tau^{2-2g} -\frac{B_{2g-2}}{2g-2}
& \text{if} \ \ g \geq 2.
\end{cases}
\end{gather}
Theorem \ref{thm.lens} follows easily from equations \eqref{eq.Ftemp7}
and \eqref{eq.Ftildeg}.

It remains to prove \eqref{eq.Ftildeg}. Equation \eqref{eq.Fglens}
and \eqref{eq.L1} implies that for $g=1$ we have
\begin{gather}
\label{eq.Ftemp4}
F_1(\tau)=f(\tau)=-\Li_1(e^\tau)-\frac{\tau}{2}-\log(-\tau).
\end{gather}
From equation \eqref{eq.Fglens}, it is easy to see that for $g \geq 1$ we have
\begin{gather}
\label{eq.Ftemp5}
F_g(\tau)= \pt^{2g-2}_{\tau} F_1(\tau)-(2g-2)! \, b_{2g-2}.
\end{gather}
Since for all $\a \in \BR$ we have
\[
\pt_\tau \Li_{\a}(e^\tau)=\Li_{\a-1}(e^\tau).
\]
Equations \eqref{eq.Ftemp4} and \eqref{eq.Ftemp5} imply that
for $g \geq 2$ we have:
\begin{gather}
\label{eq.Ftemp}
F_g(\tau)=-\Li_{3-2g}(e^\tau) +(2g-3)! \tau^{2-2g}
-\frac{B_{2g-2}}{2g-2},
\end{gather}
where we have used the fact that $(2g-2)!b_{2g}=B_{2g}/(2g-2)$.

When $g=0$, equation \eqref{eq.Fglens} implies that
\[
F_0(\tau)=\pt^{-2}_{\tau} f(\tau)
\]
up to a linear function of $\tau$, where $\pt^{-1}_{\tau}$ denotes integration
with respect to $\tau$. Integrating~\eqref{eq.Ftemp} twice using
\eqref{eq.derLi}, and matching the f\/irst three coef\/f\/icients of the
Taylor series of both sides at $\tau=0$ implies~\eqref{eq.Ftildeg} for $g=0$.
This concludes the proof of Theorem~\ref{thm.lens}.

An alternative proof of equation \eqref{eq.Ftildeg}
follows from Gopakumar--Vafa \cite[Appendix]{GV}. See also~\cite{Mr1,Mr2}.

Yet another proof of equation \eqref{eq.Ftildeg} follows from
\cite[Proposition~1]{CLZ}.\qed

\begin{remark}
\label{rem.GV}
The reader may compare Theorem \ref{thm.lens} with the BPS formulation
of~\cite{GV} and~\cite{Mr2}.
\end{remark}

\section{Future directions}
\label{sec.future}

\subsection{Analyticity of the free energy of a matrix model}
\label{sub.matrixmodel}

The notion of free energy exists for $U(N)$ gauge theories where $N$
is arbitrary. Toy models of those theories are the so-called
{\em matrix models} studied by the French school~\cite{BIPZ,BIZ}.
In a future publication \cite{GM2}, we will study the analytic properties of
the free energy of a matrix/multi-matrix model with arbitrary potential.
On the other hand, it was shown in \cite{GM1} that given a~closed 3-manifold
$M$ presented by surgery on a framed link $L$ in $S^3$, there exists
a multi-trace potential~$V_L$ whose free energy coincides with the free
energy of $M$. This observation gives geometric examples of dif\/ferent
potentials with equal free energies. Combined with the future work of~\cite{GM2}, this may give another proof of our main Theorem~\ref{thm.1}.

\subsection{Analytic continuation of the free energy in the complex plane}
\label{sub.ac}

This section is motivated by the questions that were raised in \cite{Ga2}.
Fix a closed 3-manifold $M$, and let $\rho_{M,g}$ denote the radius of
convergence of the power series $F_{\M,g}(\tau)$ at $\tau=0$. Theorem~\ref{thm.1} implies that $\inf \{ \rho_{M,g} \, | \, g \geq 0\}>0$.

\begin{question}
\label{que.1}
Fix a closed $3$-manifold $M$.
Is it true that the radius of convergence of $F_{M,g}(\tau)$ at $\tau=0$
is independent of $g$? Do the
power series $F_{M,g}(\tau)$ admit analytic continuation as multivalued
analytic functions in the complex plane minus a set of singularities,
independent of~$g$?
\end{question}

Notice that if we replace the closed 3-manifold $M$ by a knot $K$,
and consider its free energy from the loop expansion of $K$, then
the answer to this question is positive.
This was discussed in Section \ref{sub.knot}. In that case,
$F_{K,n}(\tau)$ are rational functions of $e^\tau$ and have analytic
continuation in the complex plane minus the logarithms of the roots of the
Alexander polynomial of $K$.

\begin{question}
\label{que.compute}
Can you compute the free energy (or even its planar limit) for any closed
hyperbolic manifold?
\end{question}

\subsection{The double scaling limit}
\label{sub.double}

In the physics literature, it is customary to expect that $F_g(\tau)$
has an expansion around a $g$-independent singularity $\tau_0$ as follows:
\[
F_g(\tau)= (\tau-\tau_0)^{\ga g+\delta} c_g (2g)!+O((\tau-\tau_0)^{\ga g+1}).
\]
In that case, one considers the single-variable power series
\[
f(x)=\sum_{g=0}^\infty c_g (2g)! x^g,
\]
which is called the {\em double scaling limit} of $F$. In the case
of matrix models, 2-dimensional gravity and random matrices, it turns
out that $f(x)$ satisf\/ies non-linear dif\/ferential equations which
are a specialization of the KP hierarchy. For a lengthy discussion, see
\cite{DGJ}, \cite[Section~4]{Wi} and references therein.

\begin{question}
\label{que.double}
Does the double scaling limit of the free energy of a closed $3$-manifold
satisfy a~specialization of the KP hierarchy?
\end{question}

\subsection[The $O(N)$ and $Sp(N)$ theories]{The $\boldsymbol{O(N)}$ and $\boldsymbol{Sp(N)}$ theories}
\label{sub.osp}

In the present paper, we def\/ined the free energy of a closed 3-manifold
using the $U(N)$ gauge theory. We could have used the $O(N)$ or the $Sp(N)$
gauge theory. In that case, the weight systems $W_{\mathfrak{o}_N}$ and
$W_{\mathfrak{sp}_N}$ lead to trivalent graphs that embed to unoriented
surfaces, see~\cite{B-N}. Theorem~\ref{thm.1} holds in that case, assuming
an asymptotic formula for the number of rooted unoriented maps; see Remark~\ref{rem.unoriented}. A precise asymptotic expansion
for the constants~$p_g$ of Remark~\ref{rem.unoriented} will be discussed
in a separate publication~\cite{GM3}.

\appendix

\section[Asymptotics of $t_g$ and Painlev\'e I]{Asymptotics of $\boldsymbol{t_g}$ and Painlev\'e I}
\label{sec.appendix}

In this Appendix we show that the recursion relation for $t_g$ found in~\cite{BGR} is closely related to an asymptotic, formal
solution to the Painlev\'e I dif\/ferential equation. This makes possible to
ref\/ine the asymptotic behavior obtained in~\cite{BGR}.

We f\/irst recall some results from \cite{BGR}. The asymptotics of $t_g$ is
obtained from the asymptotics of two auxiliary sequences. The f\/irst one is
$f_g$, which is def\/ined by the recursion
\begin{gather}
\label{recur}
f_g= \frac{\sqrt{6}}{96} (5g-4)(5g-6) f_{g-1} + 6 {\sqrt{6}}
\sum_{h=1}^{g-1} f_h f_{g-h},
\end{gather}
and the initial condition
\begin{gather}
f_0=-\frac{\sqrt{6}}{72}.
\end{gather}
The sequence $t_g$ is related to $f_g$ by
\begin{gather}
f_g=24^{-3/2} 6^{g/2} \Gamma\left( {5g-1\over 2}\right) t_g.
\end{gather}
Moreover, \cite{BGR} introduces another sequence $u_g$, def\/ined by
\begin{gather}
u_g= f_g \left( {25 {\sqrt {6}} \over 96} \right)^{-g} {6{\sqrt{6}}
\over [1/5]_g[4/5]_{g-1}}.
\end{gather}
Using (\ref{recur}), it is shown that $u_g$ approaches a constant $K$ as
$g\rightarrow \infty$, and this leads to the main result
of \cite{BGR} concerning the asymptotics of $t_g$. The constant $K$ is in
principle only known numerically.

To start our analysis, we def\/ine yet another sequence $a_g$ by the relation
\begin{gather}
f_g= f_0 C^g a_g, \qquad C={\sqrt {6} \over 2},
\end{gather}
so that the recursion (\ref{recur}) becomes
\begin{gather}
\label{diffeq}
a_{g} = \frac{25(g-1)^2-1}{48} a_{g-1} - {1\over 2} \sum_{\ell=1}^{g-1}
a_\ell a_{g-\ell}, \qquad a_0 = 1.
\end{gather}
We now consider the formal power series
\begin{gather}
\label{zeta0}
\phi_0(z) = z^{1\over 2} \sum_{g=0}^{\infty} a_g\, z^{-5g/2}.
\end{gather}
It is then easy to see that the recursion (\ref{diffeq}) implies that
$\phi_0(z)$ satisf\/ies the following dif\/ferential equation
\begin{gather}
\label{pone}
f^2 - {1\over 6} f'' = z,
\end{gather}
which is the well-known {\em Painlev\'e I} equation.
Therefore, the dif\/ferential equation presented at the end of \cite{BGR} is
Painlev\'e I in disguise.

We can now use properties of Painlev\'e I to give a precise description of
the
asymptotics of~$a_g$, and in turn of $f_g$ and $t_g$. This follows from the
general results of \cite{CK} together with the computations of~\cite{FIK}
for Painlev\'e I; see also \cite[Section~3]{Ka2}. We will then
content ourselves with a~statement of the asymptotic behavior of~$a_g$ and
a sketch of the main ideas.

The coef\/f\/icient $a_g$ has the following asymptotic behavior as
$g\rightarrow \infty$:
\begin{gather}
\label{aas}
a_g \sim {A^{-2g+{1\over 2}} \over \pi}\, \Gamma\left(2g-{1\over 2} \right)
S \left\{1 + \sum_{l=1}^{\infty} {\mu_{l} A^{l} \over \prod_{k=1}^{l}
(2g-1/2 -k)} \right\}.
\end{gather}
In this expression,
\begin{gather}
A ={8 {\sqrt{3}} \over 5}, \qquad S = -{3^{1\over 4} \over 2{\sqrt{\pi}}},
\end{gather}
and the $\mu_{l}$ are def\/ined by the recursion relation
\begin{gather}
\label{mul}
\mu_{l}= {5\over 16 {\sqrt {3}} l} \left\{ {192 \over 25}
\sum_{k=0}^{l-1} \mu_k a_{(l -k+1)/2} -\left( l-{9\over 10} \right)
\left( l-{1\over 10} \right)
\mu_{l-1} \right\}, \qquad \mu_0=1.
\end{gather}
It is understood here that $a_{n/2}=0$ if $n$ is not an even integer,
otherwise the coef\/f\/icient $a_{n/2}$ is given by the recursion (\ref{diffeq}).
The quantities $S$ and $\mu_{l}$ have a very clear interpretation in terms
of the dif\/ferential equation (\ref{pone}).  It is known (see for example
\cite{FIK}) that the Painlev\'e I equation exhibits a non-linear version of
the Stokes phenomenon across the line ${\rm Arg}\, z=0$.
$S$ turns out to be the corresponding Stokes parameter. The exact value of
$S$ was f\/irst obtained in \cite{Ka1} using the Riemann--Hilbert approach,
further justif\/ied by \cite{FIKN}. A calculation of the Stokes parameter using
a Borel summation approach and WKB was given by \cite{Tk}. The coef\/f\/icients
$\mu_{l}$ arise as follows. The Painlev\'e I equation admits a so-called
formal {\it trans-series solution} of the form
\begin{gather}
\phi (z) =\phi_0(z) + z^{1\over 2} \sum_{\ell=1}^{\infty}  \xi^{\ell}
\phi_{\ell}(z), \qquad \xi= C z ^{-5/8} \re^{-A z^{5/4}},
\end{gather}
where $\phi_0(z)$ is the formal power series (\ref{zeta0}), $C$ is a
constant and $\phi_{\ell}(z)$ are formal power series of the form
\begin{gather}
\phi_{\ell}(z) =\sum_{k=0}^{\infty} \phi_{\ell,k} z^{-5 k/4}.
\end{gather}
normalized by $\phi_{1}(0)=1$. The Painlev\'e I equation gives
a series of recursion relations for the coef\/f\/icients
$\phi_{\ell,k}$ in terms of $a_n$ and $\phi_{\ell',k'}$ with
$\ell'<\ell$. The coef\/f\/icients of the f\/irst series in the trans-series,
$\phi_{1,k}$, satisfy precisely the
recursion (\ref{mul}), therefore $\mu_l=\phi_{1,l}$. One f\/inds, for
the very f\/irst terms,
\begin{gather}
\phi_{1}(z)=1 - {5 \over 64 {\sqrt {3}}} z^{-{5\over 4}}
+ {75 \over 8192} z^{-{5\over 2}} - {341329 \over 23592960 {\sqrt{3}}}
z^{-{15\over 4}} +\cdots.
\end{gather}
The fact that the value of the Stokes parameter $S$ together with the
coef\/f\/icients of the f\/irst trans-series correction $\phi_{1}(z)$
lead to an asymptotic behavior of the form (\ref{aas}) has been known for a
long time in the physics literature on large-order behavior in perturbation
theory; see for example \cite{LGZJ}. A rigorous proof for the case of
Painlev\'e I was given in \cite[Section~3]{Ka2}. For a general approach to this
problem for a class of generic non-linear ODEs, see~\cite{CK}.
Using these results we can now revisit the asymptotics of the sequences
$u_g$ and $t_g$. One f\/inds indeed that
\begin{gather}
u_g \rightarrow K,
\end{gather}
as $g\rightarrow \infty$. The constant $K$ can now be determined exactly
from the asymptotics of $a_g$. Its value is,
\begin{gather}
\label{eq.Kexact}
K= {\sqrt{3\over 5}} {\Gamma(1/5) \Gamma(4/5) \over 4 \pi^2} =
0.10486898772254091800 \dots.
\end{gather}
In fact, the numerical approximation to the value of $K$ presented in
\cite{BGR} is not very precise, and can be substantially improved by
using standard techniques like
Richardson transforms. Of course, the results above give all the subleading
corrections to the asymptotics of $u_g$ (and $t_g$) in powers of $1/g$.

Let us use the Painlev\'e equation to answer a question posed by \cite{BGR}
concerning the sequence~$t_g$. The question asks whether the sequence~$(t_g)$
is {\em holonomic}, i.e., satisf\/ies a linear recursion in $g$
with coef\/f\/icients polynomials in $g$. It is easy to show that if~$t_g$ is
holonomic, then the generating series $\phi_0(z)$ satisf\/ies a linear ODE
with coef\/f\/icients polynomials in $z$. It follows that $\phi_0(z)$ (or more
precisely, its Borel transform in a suitable ray) has
analytic continuation
as a multivalued analytic function in $\BC$ minus a {\em finite} set of
singularities. On the other hand, it is known that any (truncated)
solution of Painlev\'e~I
is meromorphic (see for example~\cite{Bou,JKr}), with {\em infinitely} many
singularities that are asymptotic to one of the rays $\arg(z)=2 \pi k/5$ for
$k=1,\dots,5$ (see for example,~\cite{Bou,FIKN}).

\begin{corollary}
\label{cor.nolinear}
The sequence $(t_g)$ is not holonomic.
\end{corollary}

We end this Appendix by pointing out some connections between the results
of \cite{BGR,GJ} and the physics literature on matrix models and
two-dimensional gravity. It has been known since the pioneering work of
\cite{BIPZ,BIZ} that matrix models are a useful tool to count maps of genus~$g$
and with~$n$ vertices. The resulting generating functions at f\/ixed genus~$g$
are all analytic functions in the vertex-counting parameter, and they have
a f\/inite radius of convergence which is common to all of them. The sequence
$f_g$ in \cite{BGR} is given indeed by the coef\/f\/icients of the leading
singularities for one of these generating functions.
It is known however from the physics literature on the so-called ``double
scaling limit'' that these coef\/f\/icients are governed by the Painlev\'e~I
equation (see~\cite{DGJ}), and in fact,
as we have seen, this is the dif\/ferential equation governing the sequence~$f_g$. A study of the asymptotics of various generating functions appearing
in matrix models, including the Painlev\'e I equation, can be found in~\cite{MSW}.

\subsection*{Acknowledgements}
Much of the paper was conceived during conversations of the f\/irst and third
authors in Geneva in the spring of 2008. S.G. wishes to thank M.M. and
R.~Kashaev for the wonderful hospitality, E.~Witten who suggested that
we look at the $U(N)$ Chern--Simons theory for arbitrary $N$ and A.~Its for
enlightening conversations on the Riemann--Hilbert problem.

\newpage

\pdfbookmark[1]{References}{ref}
\LastPageEnding


\begin{thebibliography}{99}
\bibitem{Tk} Takei Y.,
        On the connection formula for the f\/irst Painlev\'e
        equation -- from the viewpoint of the exact WKB analysis,
        {\em Surikaisekikenkyusho Kokyuroku} {\bf 931} (1995), 70--99.
\end{thebibliography}
\end{document}